\documentclass[12pt, a4paper]{article}
\usepackage{amssymb,amsmath}
\usepackage[pdftex]{color}
\usepackage{pdfsync}

\newtheorem{lemma}{Lemma}

\newcommand{\tx}[1]{\mbox{\;{#1}\;}}

\numberwithin{equation}{section}
\numberwithin{lemma}{section}
\numberwithin{theorem}{section}
\numberwithin{definition}{section}
\numberwithin{remark}{section}
\numberwithin{corollary}{section}
\numberwithin{proposition}{section}
\begin{document}
\date{}
\title{Positive solutions of  semilinear elliptic problems with a Hardy potential.}
\author{Catherine Bandle\thanks{Departement Mathematik und Informatik, Universit\"at
Basel, Spiegelgasse 1, CH-4051 Basel, Switzerland,  \texttt{catherine.bandle@unibas.ch}}\:
and
Maria Assunta Pozio\thanks{Dipartimento di Matematica,  Sapienza Universit\`{a} di Roma,
P.le  A. Moro 5, I-00185 Roma, Italy, \texttt{pozio@mat.uniroma1.it}}}
\pagestyle{myheadings}
\markright{\sc }
\maketitle

\begin{abstract} Let $\Omega \subset \mathbb{R}^N$ be a bounded domain and $\delta(x)$ be the distance of a point $x\in \Omega$ to the boundary. We study the positive solutions  of the problem
$\Delta u +\frac{\mu}{\delta(x)^2}u=u^p$ in $\Omega$, where $p>0, \,p\ne 1$ and $\mu \in \mathbb{R},\,\mu\ne 0$ 
 is smaller then the Hardy constant. The interplay between the singular potential and the nonlinearity leads to interesting structures of the solution sets.  In this paper we first give the complete picture of the radial solutions in balls. In particular we establish for  $p>1$  the existence of a unique large solution behaving like $\delta^{- \frac2{p-1}}$ at the boundary. In general domains we extend results of \cite{BaPo2015} and
  show that there exists a unique singular solutions  $u$ such that $u/\delta^{\beta_-}\to c$ on the boundary for an arbitrary positive function $c \in C^{2+\gamma}(\partial\Omega) \, (\gamma \in (0,1)), c \ge 0$. Here $\beta_-$ is the smaller root of $\beta(\beta-1)+\mu=0$. 
 \end{abstract}

\bigskip

 \smallskip


\section{Introduction}
In this paper we study {\sl positive} solutions of problems of the form
\begin{align}\label{original}
 L_\mu u = \Delta u +  \frac{\mu}{\delta(x)^2}u=u^p\quad \tx{in} \Omega,
 \end{align}
where $\mu \in \mathbb{R} \setminus \{0\}$, $\delta(x)$ is the distance of a point $x\in \Omega$ to the boundary,  $p\ne 1$ is a positive constant and $\Omega \subset \mathbb{R}^N\,, \,N\ge 1$, is a bounded, connected domain with a boundary  $\partial \Omega \in C^{2+ \gamma}$, $\gamma \in (0,1)$. The expression
$$
\frac{\mu}{\delta(x)^2}=:V_\mu(x)
$$
is called the {\sl Hardy potential}.  In view of its singularities the boundary data cannot be prescribed arbitrarily. The boundary behavior depends on the interplay between the linear regime $L_\mu h=0$ and the nonlinear regime $\Delta U=U^p$. An important ingredient in the study of problem \eqref{original} is the 
{\sl Hardy constant}
\begin{align}\label{Hardy}
C_H(\Omega)= \inf_{\phi\in W^{1,2}_0(\Omega)} \frac{\int_\Omega |\nabla \phi|^2\:dx}{\int_\Omega\delta^{-2}(x)\phi^2\:dx}.
\end{align}
It is well-known \cite{MaMiPi98} that $0<C_H(\Omega)\leq 1/4$. This implies in particular that \eqref{original} cannot have nontrivial solutions belonging to $W^{1,2}_0(\Omega)$ if $\mu<C_H(\Omega)$.

A function  $h$ will be called {\sl $L_\mu$-harmonic} or simply {\sl harmonic} if it satis\-fies  $L_\mu h=0$, {\sl  sub-harmonic or super-harmonic} if $L_\mu h\geq 0$, or $L_\mu h \leq 0$, respectively. In this paper we shall only be concerned with positive sub- and super-harmonics.  

It was shown in \cite{BaMoRe09} that for $\mu\leq 1/4$, a local sub-harmonic either dominates every local super-harmonic multiplied by a suitable positive constant, or it is dominated by a multiple of any super-harmonic. This property will be referred to as the {\sl Phragmen-Lindel\"of alternative}. It was used in \cite{BaMoRe09} to determine the behavior of the solutions of \eqref{original} near the boundary. The singularity of the Hardy potential forces a solution either to vanish or to explode on the boundary. 

In the case $p>1$ the particular feature of the nonlinear problem is the existence of a {\it maximal} solution which blows up at the boundary, and in the case $p<1$ the appearance of {\it dead cores}, {\sl i.e.} regions where the solution vanishes identically.
\medskip

The structure of the radial solutions in balls is now well-understood. It has been studied in \cite{BaPo2015} in the case of sublinear nonlinearities. In order to describe the result, we have to introduce some constants which will be crucial in the sequel. For $\mu < 1/4$ we set
\begin{equation}\label{betapm}
\beta_\pm = \frac{1}{2}\pm \sqrt{\frac{1}{4}-\mu}.
\end{equation}
Furthermore let
\begin{align}\label{mu}
\mu^*= \frac{2(p+1)}{(p-1)^2}.
\end{align}
It turns out that there are essentially two types of positive solutions,
those governed by the linear regime and those with a dead core or blowup caused by the nonlinearity. More precisely we have
\bigskip

%%%%%%%%%%%%%%%%%%%TTTTTTTTTTTT
{\bf Theorem A.} {\sl  Assume $0<p<1$, $\mu<1/4, \,\mu \not= 0$ and $\Omega =B_R:=\{x\in \mathbb R^N:|x|<R\}$. 

\noindent (i) Problem \eqref{original} has a unique radial solution $u(r)$  for any $u(0)>0$. 

\noindent (ii)  For any $\rho \in (0,R)$ there exists a unique radial solution of problem \eqref{original} such that 
$$
u(r)=
\begin{cases}
0,\, \tx{if} r\in (0,\rho],\\
> 0,\, \tx{if} r\in (\rho,R).
\end{cases}
$$

\noindent (iii) There exists a unique radial solution which vanishes at $r=0$ and which is positive in $(0,R)$. Near the origin it is of the form
$$
u(r)=r^{\frac{2}{1-p}}(c" + w(r)),\quad w(0)=0 \tx{and} c"=\left(\mu^* + \frac{2(N-1)}{1-p}\right)^{1/(p-1)}.
$$
\noindent (iv)  All radial solutions satisfy $\frac{u(r)}{(R-r)^{\beta_-}}\to c>0$ as $r\to R$ and vice versa, for any constant  $c > 0$ there exists a unique solution satisfying this condition and it is radial.}
%%%%%%%%%%%%%%%%%%%%%%%tttttttttttt
\bigskip

\noindent For the superlinear case the situation is slightly different, namely
\bigskip
%%%%%%%%%%%%%%%%%%%%%%%%%%%%%%%%%%TTTTTTTTTT

{\bf Theorem B.}  {\sl Assume $p>1$, $-\mu^*<\mu<1/4, \,\mu \not= 0$ and $\Omega =B_R$.  

\noindent (i) There exists a positive number $u^*$ such that problem \eqref{original} has a unique radial solution for $u(0)\leq u^*$ whereas for $u(0)>u^*$ the radial solutions blow up before $r=R$.

\noindent (ii) If $u(0)<u^*$, $\lim_{r \to R} \frac{ u(r)}{(R-r)^{\beta_-}}=c>0$ and vice versa, for any $c > 0$ there exists a unique radial solution satisfying this condition.

\noindent (iii) If $u(0)=u^*$, then

\noindent
 $\lim_{r \to R} \frac{ u(r)}{(R-r)^{-\frac{2}{p-1}}}=(\mu+\mu^*)^{\frac{1}{p-1}}$\quad {\it  (maximal singular solution).}}

\bigskip

\noindent  In short if $u(0)<u^*$ the linear regime prevails, otherwise the nonlinearity dominates.

In the previous papers  \cite{BaPo2015,BaMoRe08} the authors considered also {\it general domains}  and constructed solutions which behaved like $c_0 \delta^{\beta_-} \leq u(x)\leq c_1 \delta^{\beta_-}$  near the boundary. In this paper we prove {\it existence and uniqueness} of solutions with a {\it prescribed boundary behavior}. The existence results in Theorem A $(iv)$ and Theorem B $(ii)$ are special cases of the following more general result
\bigskip

{\bf Theorem C.} {\sl Let $p > 0, \, p\neq 1 $, $\mu < C_H(\Omega) $ and $ \mu \not= 0$.  If $p>1$ we assume in addition that $\mu >-\mu^*$. Let $\partial \Omega \in C^{2+\gamma}$ for some $\gamma \in (0,1)$. For any   $c \in C^{2+\gamma}(\partial \Omega)$, $c \ge 0$, 
there exists a unique solution of \eqref{original} such that
\begin{align}\label{c(x)}
\lim_{\delta(x)\to 0} \Big(\frac{u(x)}{\delta(x)^{\beta_-}} -c(x^*)\Big) =0\,,
\end{align}
where $x^* \in \partial\Omega$ is the projection of $x$ on the boundary.}
\bigskip

 It was shown in \cite{BaMoRe08} that for $p>1$ and $\mu < -\mu^*$ problem \eqref{original} has no non trivial solutions. Related results are found in \cite{ DuWei2015, MaNG2017,BaPo2015, BaMaMo}.

%%%%%%%%%%%%%%%%%%%%%%%%%%%%%%%%%%%%%%%%
%%%%%%%%%%%%%%%%%%%%%%%%%%%%%%%%SSSSSSSSS
\section{Preliminaries}\label{Prelimin}
We recall some lemmas which will be used in the proofs of our theorems. 
\begin{lemma}\label{max}{\it (Maximum principle.)} {\sl  Let $\mu< C_H(\Omega)$ and $\omega\subseteq \Omega$. If $\Delta u+V_\mu u\geq 0$ in $\omega$ and $u\in W^{1,2}_0(\omega)$
 then $u\leq 0$ in $\omega$.} % the proof follows in
\end{lemma}
The proof simply follows
from the definition of the Hardy constant  in \eqref{Hardy}.
\smallskip

From this maximum principle we deduce immediately 
%%%%%%%%%%%%%%%%%%%%%%%%%%%%%%LLLLLLLLL
\begin{lemma}\label{comparison}(Comparison principle.) Let $G\subset \Omega$ be an open set such that $\overline{G}\subset \Omega$, and let $0\leq \underline{u}, \overline{u}\in W^{1,2}_{loc}(G)\cap C(G)$  be sub- and supersolutions of \eqref{original}. Assume that 
$\underline{u}\leq \overline{u}$ on $\partial G.$
\begin{itemize}
\item[(i)] If $\mu<C_H(\Omega)$, then $\underline{u}\leq \overline{u}$ in $G$.
\item[(ii)] If $p>1$ and $\overline{u}>0$ in $G$, the same statement holds without restriction on $\mu$.
\end{itemize}
\end{lemma}
%%%%%%%%%%%%%%%%%%%%%%%%%%%%%%
The fact that (ii) is valid for any $\mu\in \mathbb{R}$, was observed in \cite{BaMoRe08}. 
Also the next result is taken from \cite[Th. 2.6] {BaMoRe08}.
%%%%%%%%%%%%%%%%%PPPPPPPPPPP
\begin{lemma}\label{PrLin}{\it (Phragmen--Lindel\"of alternative)}. Let $\mu \le 1/4$ and let $\underline{h}$ be a local sub--harmonic. Then either of the following alternatives holds:
\begin{itemize}
\item[(i)] for every local super--harmonic $\overline{h} > 0$,
$$\limsup_{x\to\partial \Omega} \frac{\underline{h}}{\overline{h}} > 0\,, $$
or
\item[(ii)] for every local super--harmonic $\overline{h} > 0$,
$$\limsup_{x\to\partial \Omega} \frac{\underline{h}}{\overline{h}} < +\infty\,.$$
\end{itemize}
\end{lemma}
%%%%%%%%%%%%%%%%%%pppppppp

M. Marcus and P-T. Nguyen \cite{MaNG2017} have shown that for $0<\mu<C_H(\Omega)$ every harmonic function can be represented by the Martin kernel $K_\mu^\Omega(x,y)$, $(x,y)\in \Omega \times \partial \Omega$. For fixed $y\in \partial \Omega$, $K_\mu^\Omega(x,y)$, is an $L_\mu$-harmonic function vanishing on $\partial \Omega\setminus y$ and  equal to one at an arbitrary, but fixed point $x_0\in \Omega$. Based on estimates by Filipas, Moschini and Tertikas, Marcus and Nguyen showed that there exists a constant $c_K>1$ such that $\forall x\in \Omega, y\in \partial \Omega$,
\begin{align}\label{Martin}
c_K^{-1}\delta^{\beta_+}(x)|x-y|^{2\beta_--N}\leq K^\Omega_\mu(x,y)\leq c_K \delta^{\beta_+}(x)|x-y|^{2\beta_--N}.  \end{align}

A further important tool for establishing the existence of solutions is (cf. \cite{BaMoRe08})
%%%%%%%%%%%%%%%%%LLLLLLLLLLLL
\begin{lemma}\label{subsupexist} Let $\mu <C_H(\Omega)$ and $p\neq 1$. If there exist a sub and a supersolution 
$0\leq \underline{u} \leq \overline{u}$ in $\Omega$, then problem \eqref{original} admits a solution $U$ in $\Omega$ such that $\underline{u}\leq U \leq \overline{u}$. If $p>1$, the condition $\mu<C_H(\Omega)$ can be replaced by $\mu \leq 1/4$.
\end{lemma}
%%%%%%%%%%%%%%%
The parallel set $\Omega_\rho:=\{x\in \Omega: \delta(x)<\rho\}$ will be used to determine the behavior of the solutions near the boundary. 
If $\Omega$ is of class $C^k$, $k\geq 2$, then $\delta$ is in $C^k(\Omega_{\rho_0})$ for $\rho_0>0$ sufficiently small \cite{Foote84}. Denote as before by $x^*$ the nearest point to $x$ on $\partial \Omega$. Let $K_i(x^*)$, $i=1,..,N-1$  be the principal curvatures of $\partial \Omega$ at $x^*$.  For any $x \in \Omega_{\rho_0}$  we have
\begin{equation}\label{Deltadelta}\begin{split}
|\nabla \delta (x)| =1\,,\hspace{90pt}\\
 -\frac{N-1}{\rho_0 - \delta(x)}\le \Delta \delta(x) = - \sum_{i=1}^{N-1}\frac{K_i}{1-K_i\delta(x)} \le \frac{N-1}{\rho_0 + \delta(x)}\,.
 \end{split}
\end{equation}
The maximal distance $\rho_0$ for which $\partial \Omega_{\rho_0}$ is in $C^2$ can be estimated by means of the principal curvatures
\begin{equation}\label{rho0}\begin{split}
 0<\rho_0 \le \frac{1}{K_{max}}, \quad K_{max}:=\max\{|K_i(x)|,\,i= 1,\,\ldots\,, N-1, \,x \in \partial \Omega\}\,.
\end{split}
\end{equation}
%%%%%%%%%%%%%%%%%%%%%%%%%%%%%%%%%%%%%%%%%%%%%%%%%%%
%%%%%%%%%%%%%%%%%%%%%%%%%%%%%%%%%%%%%%%SSSSSSSSSSSSSSSS
\section{Ball}
Theorem A $(i),\,(ii),\,(iii)$ have been proved in \cite{BaPo2015}. The existence result $(iv)$ is a consequence of Theorem C.  The proof of Theorem B is based on ode techniques partly developed in \cite{BaPo2015}.
\subsection{Proof of Theorem B} 

The radial solutions of \eqref{original} in $B_R$ satisfy the ordinary differential equation ($ ':= \frac{d}{d r}$)
\begin{align}\label{radial}
u'' +\frac{(N-1)}{r} u'
+\frac{\mu}{(R-r)^2}u = u^p \quad\tx{in} (0,R), \,\,\, u'(0)=0\,.
\end{align}
This equation has for given $u(0)>0$ a unique local solution. It can be continued until one of the following cases occurs: 
\newline
1. the solution vanishes before $r=R$,
\newline 
2. it blows up before $r=R$, 
\newline 3. it exists and it is positive in the whole interval $[0,R)$. 
\smallskip

The first case is excluded by the fact that $\mu < \frac14=C_H(B_R)$. It is well-known \cite{BaMa95} that for any $0<\rho<R$ there exists a unique large solution $U_\rho$ of \eqref{radial} in $B_\rho$ which blows up at the boundary at the rate $(\mu^*)^{\frac1{p-1}}(\rho-r)^{-\frac{2}{p-1}}$. Clearly by Lemma \ref{comparison}, $U_\rho(0)$ decreases as $\rho$ increases. Define $u^*:= \lim_{\rho\to R} U_\rho(0)$. Obviously $u(r)$ blows up at $\rho<R$ if $u(0)>u^*$. This establishes the nonexistence part in Theorem B $(i)$. 
\smallskip

Next we consider the solution $U^*(r)$ of \eqref{radial} satisfying $U^*(0)=u^*$. We fix  $r_0 \in (0,R)$ and by means of suitable super and subsolutions we shall show that there exist  positive constants $0<c_1<c_2$ depending only on $r_0$, such that   $c_1(R-r)^{-\frac{2}{p-1}}\leq U^*(r) \leq c_2(R-r)^{-\frac{2}{p-1}}$ for $r\geq r_0$. 
\medskip

We start with the supersolution. Let  $w_+=c_+(R-r)^{-\frac{2}{p-1}}$. It satisfies
\begin{align}\label{eq:w1}
L_\mu w_+-w_+^p = \left[ \mu^* +\mu -c_+^{p-1} +\frac{2(N-1)}{(p-1) r}(R-r)\right]c_+(R-r)^{-\frac{2p}{p-1}}.
\end{align}
If 
$$
c_+^{p-1}>\mu^*+\mu +\frac{2(N-1)}{(p-1) r_0}(R-r_0),
$$
then $w_+$ is a supersolution of \eqref{radial} in $(r_0,R)$. 

Let $U_{r_1}(r)$ be the solution which blows up at $r_1\in (r_0,R)$.  Choose $c_+$ so large that $w_+(r_0) > U_{r_1}(r_0)$.   Let $\tilde r$ be the largest number such that $w_+(\tilde r)=U_{r_1}(\tilde r)$. Then the function 
$$
\overline u(r)= \begin{cases}
U_{r_1}(r)  &\tx{in} [0, \tilde r],\\
w_+ (r) &\tx{in} [\tilde r,R)
\end{cases}
$$
is a weak supersolution with $\overline{u}(0)>u^*$.
\medskip

Next we construct a subsolution. From \eqref{eq:w1} it follows that for small $c_-\leq (\mu+\mu^*)^{1/(p-1)}$, $w_-=c_-(R-r)^ {-\frac{2}{p-1}}$ is a subsolution in $(0,R)$.  Here we need the assumption $-\mu^*<\mu$.
Since $w_-'(0)=c_- \frac2{p-1}R^ {-\frac{2}{p-1}-1}>0$, $w_-$ is a weak subsolution as a radial solution in $B_R$.\par
\noindent
  Applying Lemma \ref{subsupexist} we deduce that there exists a solution  $\underline{u}\leq U^*\leq \overline{u}$ in $[0,R)$. The assertion $(iii)$ of Theorem B will be established if we prove the following\\
{\it {\underline{Claim}}: The function
\begin{equation*}\label{w^*}
w^*(r) := U^*(r)/(R-r)^{-\frac{2}{p-1}}
\end{equation*}
that  is  bounded from below and above by positive constants satisfies 
\begin{equation*}\label{w^*(R)} w^*(R)= \lim_{r \to R^-} w^*(r) = (\mu + \mu^*)^{\frac1{p-1}}
\end{equation*}}
In $[0,R)$,  $w^*$  satisfies the equation 
\begin{align}\label{eq:w}
(w^*_r \sigma)_r =\frac{\sigma}{(R-r)^2}&\left[ {w^*}^p-\left(\mu^*+\mu + \frac{2(N-1)}{(p-1)r}(R-r)\right)w^*\right],\\
\nonumber &\tx{where} \sigma = (R-r)^{\frac{4}{p-1}}r^{N-1}.
\end{align}
Put
\begin{align*}
f(w,r)&:=  w^p-[\mu^*+\mu+\eta(r)]w , \:w>0, \: r\in (0,R],\\
 &\eta(r)=\frac{2(N-1)}{(p-1)r}(R-r).
\end{align*} 
For any $r \in [0, R]$, the equation $f(w,r)=0$  has two zeros, namely $w=0$ and $w_0(r)= (\mu^*+\mu +\eta(r))^{\frac{1}{p-1}}$. The function $w_0(r)$ is monotone decreasing with $\lim_{r\to 0} w_0(r)= +\infty$ and $\lim _{r\to R} w_0(r) = (\mu^*+\mu )^{\frac{1}{p-1}}$.
It is important to point out that all the local maxima of the solutions of \eqref{eq:w} are below $w_0(r)$ and the local minima are above $w_0(r)$.

The  radial symmetry implies that $U_r^*(0) = 0$ and consequently 
$$
w_r^*(0)= -\frac2{p-1}U^*(0) R^{\frac{3-p}{p-1}} < 0\,.
$$
Thus $w^*$ decreases in a neighborhood of $r=0$. If it has a local minimum in $r=r_0$, by the previous remark $w^*(r_0)\geq w_0(r_0)$. 
Then $w^*(r)$ cannot have a first local maximum $w^*(r_1)$ for some  $r_1 \in (r_0 ,R)$, since we get $w_0(r_1) \ge  w^*(r_1) > w^*(r_0) \ge w_0(r_0)$, which is impossible because $w_0(r)$ is decreasing. Therefore $w^*(r)$ is either decreasing in $(0,R)$ or it has one local minimum.

Thus $w^*(r)$  is monotone in a neighborhood of $r=R$ and its bounded\-ness implies that $\lim_{r\to R} w^*(r)=w^*(R)$. Integration of \eqref{eq:w} leads to
\begin{align}\label{eq:wb}
w^*_r(r)\sigma(r)-w^*_r(r_0)\sigma(r_0)=\int_{r_0}^r\frac{\sigma}{(R-s)^2}f(w^*,s)\:ds.
\end{align}
We now distinguish between three cases.
\medskip

\noindent 1. If $p<5$ we can pass to the limit as $r\to R$ and obtain $\lim_{r\to R} w^*_r\sigma =\ell$. Since $w^*$ is bounded, $\ell=0$ (otherwise $\sigma ^{-1}$ would not be integrable in a neightbourhood of $R$). Hence
$$
-w_r^*(r)\sigma(r)=\int_r^R\frac{\sigma}{(R-s)^2}f(w^*,s)\:ds.
$$
Elementary calculus shows that $w^*$ is uniformly bounded only if\newline
 $f(w^*(R),R)= 0$.  
 \medskip

\noindent 2.  If $p = 5$ then $\sigma (r) = (R-r)r^{N-1}$. It was observed before that $w^*(r)$ is bounded and monotone near $r=R$. If in such a neighborhood $w_r\geq 0$ and  $\liminf_{r\to R} w_r^*(r)\sigma(r)>0$,   $w$ would be unbounded since $\sigma^{-1}$ is not integrable near $r=R$. {\sl i.e.} a contradiction. Analogously in case $w_r\leq 0$ we get a contradiction if $\limsup_{r\to R} w_r^*(r)\sigma(r)<0$. Thus in  case $w_r\geq 0$ we have $\limsup_{r\to R} w_r^*(r)\sigma(r)\ge 0$ and $\liminf_{r\to R} w_r^*(r)\sigma(r)=0$ . In case $w_r\leq 0$, we  have  $\limsup_{r\to R} w_r^*(r)\sigma(r)=0$ and $\liminf_{r\to R} w_r^*(r)\sigma(r) \le 0$.  Moreover these results together with \eqref{eq:wb} show that there exists a sequence $r_n \nearrow R$ as $n \to \infty$ such that
$$
\int_{r_0}^{r_n}\frac{\sigma}{(R-s)^2}f(w^*,s)\:ds = \int_{r_0}^{r_n}\frac{s^{N-1}}{(R-s)}f(w^*,s)\:ds 
$$
has a finite limit for $n \to \infty$. Hence  $f(w^*(R),R)= 0$.
\medskip

\noindent 3. If $p > 5$, we divide \eqref{eq:wb} by $\sigma$ and integrate. We then obtain
\begin{gather*}
w^*(r) - w^*(r_0) - w_r^*(r_0)\sigma(r_0)\int_{r_0}^r\sigma^{-1}(s)\:ds\\
 = \int_{r_0}^r\frac{\sigma(s)}{(R-s)^2}f(w^*,s) \big(\int_{s}^r\sigma^{-1}(y)\:dy\big)\:ds \,.
\end{gather*}
Again it follows from elementary calculus that $w^*(R)$ is bounded only if $f(w^*(R),R)=0$. The case $w^*(R)=0$ is excluded because $w^*(r) >c_->0$. Consequently $w^*(R)=w_0(R)$ which establishes Theorem B $(iii)$. Moreover  $w^*$ decreases in $[0, R]$ because if $w^*$ has a local minimum at $r_0 \in (0, R)$, then $w(R)=  (\mu+\mu^*)^{1/(p-1)} > w(r_0) > (\mu+\mu^*)^{1/(p-1)}$, which is impossible.
\medskip

Next we proceed to the proof  Theorem B $(ii)$.
Consider the function $w(r)=u(r)(R-r)^{2/(p-1)}$  where $u(r)$ solves \eqref{radial} and $u(0)<u^*$. 

We claim that  $w<w^*$ in $(0,R)$ and that $w(R)=0$.  In fact $w$ cannot intersect $w^*$. Namely if $w(r_1)=w^*(r_1)$ and $w<w^*$ in $(0,r_1)$, then 
$U^*>u$ in $B_{r_1}$. This implies $L_\mu (U^*-u)>0$. From Lemma \ref{max} it follows that $U^*<u$, which is obviously a contradiction. This establishes the first assertion.
The same arguments as for $w^*$ show that $w< w^*$ decreases in the whole interval $(0,R)$ and that we must have $f(w(R),R)=0$. We shall prove that $w(R)<w^*(R)$, which  implies that
 $w(R)=0$. 
 
Suppose that $w(R) = w^*(R)$. Then in addition to $w^*$ there is another solution satisfying $(iii)$, namely  $u(r) = w(r)(R-r)^{-2/(p-1)}< U^*(r) = w^*(r) (R-r)^{-2/(p-1)}$. Define $m(x)$ by $u(r)= m(r) U^*(r)$. Then $m \le 1 $ in $B_R$ and $m = 1 $ on $\partial B_R$ and, since $u$ solves \eqref{original} in $B_R$, the differential equation for $m(r)$ is
\begin{align}
 U^*\Delta m + 2 <\nabla m,\nabla U^*>+ m \Delta U^* +  \frac{\mu}{\delta^2}m U^*=m^p (U^*)^p.
 \end{align}
Hence
\begin{align}
 U^*\Delta m + 2 <\nabla m,\nabla U^*>=m^p (U^*)^p - m (U^*)^p \le 0.
 \end{align}
The maximum principle implies that $m$ attains its minimum on the boundary. This is only possible if $m(r) \equiv 1$, i.e. $U^*$ is the unique solution satisfying $\lim_{r\to R} U^*(r)(R-r)^{2/(p-1)}=(\mu +\mu^*)^{1/(p-1)}.$ Notice that uniqueness result holds also for non radial solutions in general domains.
\medskip
 
We now need to prove the boundary behavior in $(ii)$.\par
\noindent The function $u(r)=w(r)(R-r)^{-\frac{2}{p-1}}$ is a solution of \eqref{radial} which can be written as follows
\begin{align}\label{eq:uw}
u'' +\frac{N-1}{r}u' + \frac{\mu-w^{p-1}}{(R-r) ^2}u=0, \:u'(0)=0.
\end{align}
A radial harmonic in $B_R$ solves
$$
h''(r) +\frac{N-1}{r} h'(r) + \frac{\mu}{(R-r)^2} h(r)=0\,\, \tx{in} [0,R),\: h'(0)=0.
$$
It is a supersolution to \eqref{eq:uw} and by Lemma \ref{comparison} it is a lower bound for $u$ provided $h(0)<u(0)$. We can now apply the Fuchsian theory to $h$. Since $\mu<1/4$,  $\lim_{r\to R} h(r)/(R-r)^{\beta_-}=c_0>0$. Consequently $u(r)\geq c_1(R-r)^{\beta_-}$ for some sufficiently small constant $c_1$.

In order to construct an upper bound we take the solution $u_\epsilon$ of  \eqref{eq:uw}  with $w$ replaced by $w_\epsilon= \max\{w,\epsilon\}$.  It is a subsolution  to \eqref{eq:uw} and if $u_\epsilon(0)>u(0)$ it is an upper bound for $u$.  Since $w$ is decreasing and $w(R)=0$, we have $w_\epsilon =\epsilon$ near $r=R$. We can now apply the Fuchsian theory to  $u_\epsilon$.
The indicial equation for $u_\epsilon$ implies that $u_\epsilon(r)=(R-r)^{\beta_\epsilon} f(R-r)$ where $\beta_\epsilon$ is the smaller root of  $\beta_\epsilon(\beta_\epsilon-1) +\mu -\epsilon^{p-1}=0$ and $f$ is an analytic function near $r=R$ such that $f(0)=b>0$. Hence for some $\delta_0 >0$  we have 
\begin{equation}\label{uepsilEstimate}
0<c_0 (R-r)^{\beta_-} \le u(r) \le u_\epsilon(r)\le  k (R-r)^{\beta_\epsilon}\,, \,\, \forall  r \in (R-\delta_0, R)\,.
\end{equation}
We replace $r$ by $\delta = R-r$ and consider the function $ v(\delta) = u(R-\delta) \delta^{-\beta_-} $. A straightforward computation leads to
\begin{align}\label{Eqv} 
(\sigma  v')'= \sigma \left(v^p\delta^{\beta(p-1)}+\beta\frac{N-1}{(R-\delta)\delta}v\right),\\
 \tx{where}\,\,
 \nonumber \sigma(\delta)=\delta^{2\beta}(R-\delta)^{N-1}\,\,  \tx{and}\,\,  \beta = \beta_-. 
\end{align}

We take $0 < \delta < \delta_0$ and we integrate \eqref{Eqv}  to get:
\begin{align}\label{eq:v4}
v(\delta)-v(\delta_0) +\sigma(\delta_0) v'(\delta_0)\int_\delta^{\delta_0} \sigma^{-1}\:ds =\\
 \nonumber \int_\delta^{\delta_0}\sigma (v^ps^{\beta(p-1)}+\beta\frac{N-1}{(R-s)s}v)(\int_\delta^s\sigma^{-1}\:d\xi)\:ds.
\end{align}
 Keeping in mind that $\beta <1/2$, we have for $s \in (\delta, \delta_0]$  
$$
\int_\delta^s \sigma^{-1}\:d\xi \leq \frac{s^{1-2\beta }}{(1-2\beta)(R-\delta_0)^{N-1}}<\infty.
$$
From \eqref{uepsilEstimate} and the definition of $v$ it follows that $v \le k \delta^{\beta_\epsilon-\beta}$. Hence there exist some positive constants $k_1,\, k_2,\, k_3$ independent of $\delta$  such that
\begin{align*} \int_\delta^{\delta_0}\sigma (v^ps^{\beta(p-1)}+\beta\frac{N-1}{(R-s)s}v)(\int_\delta^s\sigma^{-1}\:d\xi)\:ds \hspace{60pt} \\ \le k_1 \int_\delta^{\delta_0} s^{2\beta + \beta_\epsilon-\beta-1+1-2\beta}\: ds 
+ k_2  \int_\delta^{\delta_0} s^{2\beta + p(\beta_\epsilon-\beta)+ \beta (p-1) +1-2\beta}\: ds \le k_3 \,,
\end{align*}
for sufficiently small $\epsilon$. By our assumption $\mu>-\mu^*$ we have $\beta>-\frac{2}{p-1}$. It is therefore possible to choose $\epsilon$ so small that $p(\beta_\epsilon-\beta)+ \beta (p-1) +1>-1$, and $\beta_\epsilon-\beta >-1$. Hence the integrals above converge as $\delta \to 0$
and $v$ is uniformly bounded.
\bigskip

Next we want to show that $v(\delta)$ has a limit as $\delta$ tends to $0$.

\noindent If $\beta = \beta_- >0$, \eqref{Eqv}  implies that $v$ cannot have a local maximum and since $v'(R) = u'(0) R^{-\beta}  -\beta u(0) \delta^{-\beta-1} = - \beta u(0) \delta^{-\beta-1} <0$, it decreases in the whole interval. Then $u(r)/(R-r)^{\beta_-}$ is an increasing bounded function for $ r = R-\delta \in (0, R)$, hence it has a limit as $r \to R^-$ .\par
\noindent Assume $\beta =\beta_- <0$. Recall that $v$ is uniformly bounded. Since by assumption $\beta_->-\frac{2}{p-1}$ all integrals in  \eqref{eq:v4} exist for $\delta=0$. Consequently we can pass to the limit which shows that $\lim_{\delta \to 0} v(\delta)$ exists.
\medskip

 For the last assertion of $(ii)$ we refer to the proof of  Theorem C. This completes the proof of Theorem B.\hfill      $\square$

%%%%%%%%%%%%%%%%%%%%%%%%%%%%%%%%%%%%%%%%%%%%%%%%%%%%
%%%%%%%%%%%%%%%%%%%%%%%%%%%%%SSSSSSSSSSS
\section{General domains}

%%%%%%%%%%%%%%%%%%%%%%%%%%%%%sssssssssssssssssss
\subsection{ Proof of Theorem C} The existence of a solution $u$ of \eqref{original} satisfying \eqref{c(x)} will be proved  by means of a supersolution 
 $\overline{u}$ and a subsolution  $\underline{u}$ of \eqref{original},  which both satisfy \eqref{c(x)} and which are such that $ \underline{u} \le \overline{u}$.\par
We start with an important observation. For  $\delta\le \delta_0 \le \rho_0/2$, $\rho_0$ defined in Section \ref{Prelimin}, set
\begin{equation}\label{deltaw}
u(x)= \delta^{\beta}w(x)\,,\qquad(\beta = \beta_-).
\end{equation}
Then $u$ is a local solution of \eqref{original} in $\Omega_{\delta_0}$, if and only if  $w$ satisfies
 \begin{equation}\label{zeta}
\begin{array}{c}
\mathcal A(w):=\Delta w +\frac{2 \beta}{\delta}<\nabla \delta, \nabla  w> + \frac{\beta \Delta\delta}{\delta } w=\delta^{\beta(p-1)}w^p\,.\\
 \end{array}
\end{equation}
The function $\underline{u}=\delta^{\beta}\underline{w}(x)$ is a local  subsolution of \eqref{original} in $\Omega_{\delta_0}$ if \eqref{zeta} holds with the equality sign replaced by  $\geq$.
Analogously $\overline{u}=\delta^\beta \overline{w}$ is a super solution if the inequality sign is reversed.\par
\medskip

\noindent We shall construct a supersolution as the minimum between a local supersolution satisfying \eqref{c(x)} and a global one which satisfies  \eqref{c(x)} with $\ge$ instead of equality.\par\noindent
{\it Local supersolution}.

Let $\alpha \in (0,1)$ be such that
\begin{equation}\label{alpha}
\alpha \in \begin{cases} 
(0,1) &\tx{if} \beta_-<0,\\
 (0,1-2\beta)&\tx{if}\beta_->0.
 \end{cases}
\end{equation}
 For any given 
$c \in C^{2+\gamma}(\partial \Omega)$,
 $c \ge 0$, let  $h \in C^{2+\gamma}(\overline\Omega)$ be the solution  %(cf. \cite{DiB} Ch.2, Th. 6.1)
of 
\begin{equation}\label{h}
\left\{
\begin{array}{c}
\Delta h = 0\,,\,\,\, x \in \Omega\,,\qquad\\
h(x) =c(x)\,,\,\, x \in \partial\Omega\,.
 \end{array}\right.
\end{equation}
Set $\overline w = h(x) + A \delta^\alpha$ where $A\ge 1$ will be fixed below. 
Then for $x \in \Omega_{\delta_0}$ 
{\begin{align*}
{\mathcal A(\overline w)} =  -A \alpha(1-\alpha) \delta^{\alpha -2} +A \alpha \delta^{\alpha -1} \Delta \delta +\frac{2 \beta}{\delta}<\nabla \delta, \nabla  h>\\
 +2A \beta\alpha\delta^{\alpha-2} + \frac{\beta \Delta\delta}{\delta}\overline w.
 \end{align*} 

Since $|\Delta \delta| \leq c$ and $|\nabla\delta|\equiv1$ in $\Omega_{\delta_0}$  we 
have
 \begin{align*}\label{sopraw}
 {\mathcal A}(\overline {w})\le
 A\alpha(\alpha +2\beta-1)\delta^{\alpha-2} +
\nonumber A(\alpha + |\beta|) c\delta^{\alpha-1}+\frac{2|\beta|}{\delta}|\nabla h|_\infty +\frac{|\beta| c}{\delta} |h|_\infty.
\end{align*}
By our assumptions on the regularity of $\partial \Omega$ and $c$, $|h|$ and $|\nabla h|$ are uniformly bounded in $\Omega$, cf. \cite[page 161]{mi}. Hence  $2|\beta| |\nabla h|_{\infty} + |\beta| c |h|_{\infty} <C$ and consequently
$$
 {\mathcal A}(\overline {w})\le \delta^{\alpha-2}A[\alpha(\alpha+2\beta-1)+ \delta (\alpha+|\beta|) c + \frac{C}{A}\delta^{1-\alpha}].
 $$
 Since $\alpha$ satisfies \eqref{alpha}, it is possible to choose $\delta_0<\rho_0$ sufficiently small such that \begin{equation}\label{soprawbis}
{\mathcal A(\overline w)} \le 0 \le \delta^{\beta(p-1)}\overline w^p\,,\,\, x \in \Omega_{\delta_0},\,\,\forall\, A>1.
\end{equation}
Thus $\hat{u}:=\delta^\beta \bar{w}$ is a local super--harmonic function and therefore the desired local supersolution.
\bigskip

Next we construct a supersolution in the whole domain. We shall treat the case $\mu>0$ and $\mu<0$ separately.\medskip

\noindent
{\it Global supersolution  for $\mu <0$.}  Let  $\eta=\eta(r)$ be the solution of
\begin{equation}\label{etaA}
\left\{\begin{array}{c}
\eta'' +\frac{(N-1)}{r}\eta' + \frac{\mu}{(\delta_0-r)^2}\eta = 0\,,\,\, r \in (\delta_0/2, \delta_0),\\
\eta(\delta_0/2)= 1\,,\,\, \eta'(\delta_0/2)=0\,.  \hspace{35pt}                                                                                    \end{array}\right.
\end{equation}
Since $\mu<0$ the function $\eta(r)$ is increasing in a neighborhood of $\delta_0/2$ and it has no local maximum. Thus $\eta(r)$ is a positive increasing solution in $(\delta_0/2, \delta_0)$. 

We claim that
\begin{equation}\label{bryUn}
\lim_{r \to \delta_0} \frac{\eta(r)}{(\delta_0-r)^{\beta_-}} = C_\eta >0 \,.
\end{equation}
For the proof of this claim we proceed as in \cite{BaPo2015} where a similar result has been derived for the nonlinear equation. We choose  $\delta=\delta_0-r$
 instead of $r$ as the new variable and
 set $\eta(\delta_0-\delta) =\delta^{\beta} v$ where $\beta=\beta_-<0$.
From \eqref{etaA} we obtain
\begin{align}\label{Ev0}
 v''+ \left (2\frac{\beta}{\delta}-\frac{N-1}{\delta_0-\delta}\right )v'-\beta\frac{N-1}{(\delta_0-\delta)\delta}v
= 0\, ,\,\, \,  \delta \in (0, \delta_0/2)\,.
\end{align}
This equation can be written in the form
\begin{align}
(\sigma v')'= \beta \sigma \frac{N-1}{(\delta_0-\delta)\delta}v , \tx{where}
 \sigma(\delta)=\delta^{2\beta}(\delta_0-\delta)^{N-1}.\label{Ev1}
\end{align}
Since $\beta<0$ it follows that  $(\sigma v')' \le 0$ which implies that for $\epsilon <\delta$
$$
\sigma(\epsilon) v'(\epsilon)  \geq \sigma(\delta) v'(\delta)\:.
$$
Dividing the last inequality by $\sigma(\epsilon)$ and integrating in the same interval we get
\begin{align*}
0<v(\epsilon)\leq v(\delta)-\sigma(\delta) v'(\delta)\int_\epsilon^{\delta}\sigma^{-1}(s)
\:ds\,.
\end{align*}
This  implies that $v$ is bounded as $\epsilon \to 0$.\par
\noindent
Integrating  \eqref{Ev1} we get
$$
v'(\epsilon) = \sigma^{-1}(\epsilon) \sigma(\delta) v'(\delta) - \beta\sigma^{-1}(\epsilon)\int_\epsilon^{\delta}\sigma(s)\frac{N-1}{(\delta_0- s)s}v\,ds \,.
$$
We can pass to the limit  $\epsilon \to 0$ and conclude, since $\beta_-$ is negative,
 that $|v'|$ is bounded. Hence there exists 
$$\lim_{\delta \to 0}v(\delta)= v(0)\,.$$
 Suppose that our claim \eqref{bryUn} is not true and that  $v(0) = 0$. Following the proof of Lemma 2.4 in \cite{BaPo2015}, we define $w$ by $\eta= \delta^{\beta_-} v = \delta^{\beta_+}w$. It satisfies the same equation \eqref{Ev1} as $v$ with  $\beta_-$ replaced by $\beta_+$,
\begin{equation}\label{Ew1}\begin{split}  
(\sigma_+ w')'= \beta_+ \sigma_+ \frac{N-1}{(\delta_0-\delta)\delta}w (\ge 0), \\
\tx{where}
 \sigma_+(\delta)=\delta^{2\beta_+}(\delta_0-\delta)^{N-1}.
\end{split}
\end{equation}
We integrate \eqref{Ew1} over  $[\epsilon, \delta]$ for some $\delta \in (0, \delta_0/2]$  and we get
\begin{align*}
\sigma_+(\delta) w'(\delta) - \sigma_+(\epsilon) w'(\epsilon) =\beta_+ \int_\epsilon^\delta \sigma_+(s) \frac{N-1}{(\delta_0-s)s}w(s)\,ds=\\
 \beta_+ (N-1) \int_\epsilon^\delta(\delta_0-s)^{N-2}v(s)\,ds\,.\label{Ew2}
\end{align*}
Since $v$ is bounded we can pass to the limit as $\epsilon \to 0$ and deduce that  $\lim_{\epsilon \to 0} \sigma_+(\epsilon) w'(\epsilon) = M$. and hence $w(\delta) = O( \delta^{1-2\beta_+})$.  By assumption $v(0)=0$ and therefore $\delta^{2\beta_+ -1} w = v \to 0$ as $\delta \to 0$, i.e. $M=0$. Thus
\begin{align*}
 w'(\delta) =\beta_+ \sigma_+^{-1}(\delta) \int_0^\delta \sigma_+(s) \frac{N-1}{(\delta_0-s)s}w(s)\,ds\geq 0\:.
\end{align*}
Hence $w(\delta) \ge w(\epsilon) \ge 0$, i.e. $w(\epsilon)$ is bounded in a neighborhood of $\delta =0$. This implies $\eta(\delta_0 -\delta) \le c \delta ^{\beta_+} \to 0$ as $\delta \to 0$.  This is impossible since $\eta$ increases as  $\delta\to 0^+$. Consequently $v(0)>0$ and  \eqref{bryUn} follows.
\bigskip

Choose
\begin{equation}\label{Mchoice}
M>\frac{|c(\cdot)|_\infty}{C_\eta},
\end{equation}
where $C_\eta$ is defined in \eqref{bryUn}.
Then the function
\begin{equation}\label{tildeu}
\tilde{u}(x): = \left\{\begin{array}{c}
            M\eta(\delta_0 -\delta(x))\,,\,\,x \in \Omega_{\delta_0/2}\,,
             \\
            M \,,\,\,x \in \Omega \setminus \Omega_{\delta_0/2}\,.
          \end{array}\right.
\end{equation}
is in $C^1(\Omega)$ and is a (weak) supersolution of \eqref{original} satisfying
\begin{equation}\label{c(x)m}
\liminf_{\delta(x, \partial \Omega)\to 0} \frac{\tilde{u}}{\delta(x)^{\beta_-}} > |c(\cdot)|_\infty.
\end{equation}
Indeed since $\mu <0$, any nonnegative constant is a supersolution. Moreover since $\eta'(r) \ge 0$ and $\delta_0 \le \rho_0/2$, by \eqref{Deltadelta} we have
$ - \eta' \Delta \delta \le \frac{(N-1)}{\delta_0-\delta} \eta'$. Hence
\par\noindent
$$\Delta\tilde{u}+ \frac{\mu}{\delta^2}\tilde{u}= M\Big(\eta''-\eta' \Delta\delta + \frac{\mu}{\delta^2}\eta\Big) \le M\Big(\eta'' +\frac{(N-1)}{\delta_0 -\delta}\eta' + \frac{\mu}{\delta^2}\eta\Big) =0 \le \tilde{u}^p\,.
$$
\bigskip

\noindent
{\it Global supersolution for $\mu \in (0, C_H(\Omega))$}. 
We consider the function $z(x) = \int_{\partial\Omega}K^\Omega_\mu(x,y)\,dS_y$, where $K^\Omega_\mu$ is the Martin Kernel introduced in Section \ref{Prelimin}.  By \cite{MaNG2017}  we have that $z$ is  harmonic in $\Omega$.  Next we will prove that estimate \eqref{Martin} implies
\begin{equation}\label{c0}
\liminf_{d(x, \partial \Omega)\to 0} \frac{z}{\delta(x)^{\beta_-}} >c_0 > 0\,.
\end{equation}
For any $x \in \Omega _{\delta_0}$, and $y \in \partial \Omega \cap B_{\delta(x)/2}(x^*(x))$, there holds
$$
 |x-y| \le |x-x^*(x)| + |x^*(x)-y| \le \delta(x) + \delta(x)/2 = 3 \delta(x)/2\,,
 $$
where $x^*(x)\in \partial \Omega$ is the nearest point to $x$.
Obviously 
$$
z (x) \ge \int_{\partial\Omega\cap  B_{\delta(x)/2}(x^*(x))}K^\Omega_\mu(x,y)\,dS_y.
$$
The boun\-dary re\-gularity together with the fact that $x \in \Omega_{\delta_0},$ $(\delta_0<\rho_0/2)$, imply that $|\partial\Omega\cap B_{\delta(x^*(x))/2}| \ge \epsilon_0 \delta(x)^{N-1}$ for some $\epsilon_0> 0$.
Then  by  \eqref{Martin} there exists a constant $c_0>0$ independent of $x$ such that
\begin{equation}
z (x) \ge c_K^{-1}\delta(x)^{\beta_+} \big( 2 \delta(x)/3\big)^{2\beta_- - N}\epsilon_0 \delta(x)^{N-1} >c_0 \delta(x)^{\beta_-}\,.
\end{equation}
For a given $M \ge \frac{|c|_\infty}{c_0}$ define
$$
 \tilde{u}(x):= M z(x).
 $$
Then $\tilde{u}(x)$ satisfies \eqref{c(x)m} and 
$$ L_\mu \tilde{u}(x) =0 \le \tilde{u}(x)^p\,.$$
Hence $\tilde{u}(x)$ is a supersolution of \eqref{original} satisfying $\tilde u\geq c(x)$ on $\partial \Omega$.
\medskip

Next we construct  supersolutions in the whole domain with boundary values $c(x)$, cf \eqref{c(x)}.
\bigskip

\noindent
{\it  Supersolutions satisfying \eqref{c(x)}}
Consider  the local supersolution $\hat{u}=\delta^\beta \overline{w}$ constructed at the beginning of Section 4.2. Choose
$$A >  (\delta_0/2)^{-\beta_--\alpha} \max\{ \tilde{u}(x)\,:\,\, x\in \Omega \,,\, \delta(x)=  \delta_0/2\},$$
 where $\tilde{u}$ is the global supersolution satisfying $\tilde u\geq c(x)$ on $\partial \Omega$.
 Then for any $x\in \Omega$ such that $\delta(x)=\delta_0/2$ we have $\tilde{u}(x) \le \delta_0^{\beta_-+\alpha}A \le \hat{u}(x)$, while $\tilde{u} >\hat{u}$ in a neighbourhood of $\partial \Omega$.  Then
\begin{equation}\label{overu}
\overline{u}(x): = \left\{\begin{array}{c}
           \min\{\hat{u}(x),\,\tilde{u}(x)\}\,,\,\,x \in \Omega_{\delta_0/2}\,,
             \\
            \tilde{u}(x) \,,\,\,x \in \Omega \setminus \Omega_{\delta_0/2}\,.
          \end{array}\right.
\end{equation}
is the desired supersolution which satisfies \eqref{c(x)}.
\bigskip

\noindent
{\it Subsolution  satisfying \eqref{c(x)}.}
Set $\underline w = (h(x) - a \delta^\alpha)_+$ where $\alpha$ is defined in \eqref{alpha},  $h$ solves \eqref{h} and $a>0$  will be fixed below so that the support of $\underline w$ is contained in $\Omega_{\delta_0}$. This is the case if $a \ge a_0 >0$ for a suitable $a_0$. If  ${\mathcal A(\underline w)} \ge  \delta^{\beta(p-1) }\underline w^p $, at points where $\underline w >0$  and since $0$ is a solution of \eqref{original} then,
\begin{equation} \underline{u} =\,\left\{\begin{array}{c} \delta^\beta \underline w\,,\, \,\,x \in \Omega_{\delta_0}\,,\\
0\,,\,\,\,x\in \Omega\setminus \Omega_{\delta_0}\,,
\end{array}\right.
\end{equation}
is the desired subsolution. Since $a \delta^\alpha \le |h|_\infty$ where $\underline w >0$, we have
\begin{align}\label{sottou}
\mathcal A(\underline w) &=  a \alpha(1-\alpha) \delta^{\alpha -2} - a \alpha \delta^{\alpha -1} \Delta \delta +\frac{2 \beta}{\delta}<\nabla \delta, \nabla  h>
-2a \beta\alpha\delta^{\alpha-2} + \frac{\beta \Delta\delta}{\delta}\underline w\\
\nonumber &\ge a \alpha (1-\alpha-2\beta)\delta^{\alpha -2}- \delta^{-1}[(\alpha + |\beta|) |h|_\infty |\Delta \delta|_\infty +2|\beta||\nabla h|_\infty ]\,.
 \end{align}
In $\Omega_{\delta_0/2}$, where $\underline{w}>0$ we want to have
$$ \delta^{2-\alpha} {\mathcal A(\underline w)} \ge \delta^{\beta(p-1) + 2-\alpha}\underline w^p\,,$$
If $p<1$ it suffices to choose $a$ sufficiently large, whereas in the superlinear case $p>1$ we have to require in addition that $\beta(p-1) + 2-\alpha >0$. In view of our assumption $\mu>-\mu^*$ this condition can always be satisfied by choosing 
$\alpha$ very small.
By construction the subsolution is below the supersolution. Hence the existence of a solution satisfying \eqref{c(x)} follows. 
\vskip0,6cm

\noindent
{\it Uniqueness}.\, For any given boundary data $c \in C^{2+\gamma}(\partial\Omega)$, $c \ge 0$, let $u_1,\,u_2$ solve \eqref{original} and \eqref{c(x)}. We argue by contradiction. Suppose that  $\underline{h}:=(u_1  - u_2)_+\not\equiv 0$ (or  $\underline{h}:=(u_2 - u_1)_+\not\equiv 0$). Clearly   $\underline{h}$}  is a sub--harmonic function, and it satisfies $\lim_{x\to\partial \Omega} \frac{\underline{h}}{\delta^{\beta_-}} =0$.\par
In \cite{BaMoRe08}  simple local super-harmonic functions $\overline H$ and $\overline{h}$ have been constructed with the property that  

\noindent $\lim_{d(x,\partial \Omega)\to 0} \frac{\overline H}{\delta^{\beta_-}} =1$ and $\lim_{d(x,\partial \Omega)\to 0} \frac{\overline h}{\delta^{\beta_+}} =1$. 

\noindent Since $\limsup_{x\to \partial \Omega} \frac{\underline{h}}{\overline{H}}=0$  the Phragmen--Lindel\"of alternative $(ii)$ (see Lemma \ref{PrLin}) applies and we get $\underline{h}<c\delta^{\beta_+}$. 
Consequently  $(u_1  - u_2)_+$ belongs to $W^{1,2}_0(\Omega)$. Since $\mu < C_H(\Omega)$ the comparison  principle implies that $(u_1  - u_2)_+ \equiv 0 $ in contradiction to our assumption.

 \noindent
 This completes the proof of Theorem C. \hfill $\square$
 %%%%%%%%%%%%%%%%%%%%
 \bigskip
 
 {\bf Open problem}
It is not known what is the behavior of the solution in Theorem C at the points where $c(x^*)=0$. According to a local version of the Phragmen-Lindel\"of alternative, which is not yet available, it should be $\lim_{\delta(x)\to 0} u(x)/\delta(x)^{\beta_+} =c$.
%%%%%%%%%%%%%%%%%%%%%%%%%%%%%%%%%%%%%%%%%%%%%%%%%%%%%%%%%%%%%%%%%%%%%%%%%%%%%%%%%

\end{document}